%Paper: 9202010
%From: JLITTLE@hcacad.holycross.edu
%Date: Thu, 13 Feb 1992 09:22 EDT

%This is a plain TeX file.  It uses no special fonts, etc.

\magnification = \magstep1
\def\bC{{\bf C}}
\def\bP{{\bf P}}
\def\bZ{{\bf Z}}
\def\eps{\varepsilon}
\def\lra{\longrightarrow}
\vglue 1truein
\centerline{\bf Another Relation Between Approaches}
\centerline{\bf to the Schottky Problem}
\bigskip
\centerline{John B. Little}
\bigskip\bigskip
\beginsection 1. Introduction \par
The recent extensive work on several approaches to the Schottky problem
has produced marked progress on several fronts.  At the same time, it has
become apparent that there exist very close connections between the various
characterizations of Jacobian varieties described in Mumford's classic
lectures [M]
and the more recent approaches related to the K.P. equation.  Some of the
most striking results of this kind are to be found in the papers [B-D] and
[F].  In the first of these, Beauville and Debarre show that for a
principally-polarized abelian variety (p.p.a.v.) $(A,\Theta)$ of
dimension $g$, the Andreotti-Mayer condition
${\rm dim}\ \Theta_{Sing} \ge g - 4$ is a consequence of any one of the
following hypotheses:
\medskip
\item{$\bullet$}  There are distinct points $x,y,z \in A$ such that
$\Theta \cap \Theta_z \subset \Theta_x \cup \Theta_y$,
\item{$\bullet$}  The Kummer variety $K(A)$ has a trisecant line, and
\item{$\bullet$} The theta function $\theta_A$ satisfies the K.P.
equation (in the sense described following the Main Theorem below).
\medskip

In the second paper cited above, Fay indicates a possible relation between
the classical Schottky-Jung relations and the K.P. hierarchy.

However, until now, the approach via double translation manifolds has
seemed to be quite different from these other approaches to the
Schottky problem.  The purpose of this paper is to bring this last
approach ``into the fold'' as it were, and to show precisely how
it relates to the approaches via trisecants and flexes of the Kummer
variety, and via the K.P. equation.  Another relation was pointed
out in [L], but we believe that the following theorem gives a much more
complete indication of the precise connection.

\proclaim Main Theorem.  Let $(A, \Theta)$ be an indecomposable
$g$-dimensional p.p.a.v. over $\bC$, and let
$\psi_A \colon A \rightarrow \bP^{2^g - 1}$
be the mapping defined by the linear system of second-order theta
functions on $A$.  Assume that the image $K(A) = \psi_A(A)$ has
a ``curve of flexes,'' or more precisely that there is an irreducible
curve $\Gamma \subset A$ such that for generic $p \in \Gamma$
$$\Gamma - p \subset V_{Y_p} = \{ 2t \in A \vert t + Y_p \subset \psi_A^{-1}
                               (\ell)\ {\rm for\ some\ line\ } \ell \subset
                               \bP^{2^g - 1}\},$$
where $Y_p$ is the artinian length-3 subscheme
${\rm Spec}({\cal O}_{\Gamma,p}/m_{\Gamma,p}^3) - p \subset (A,0)$.
Then $\Theta$ is a generalized translation manifold.  That is,
$\Theta$ has a local parametrization of the form:
$$z_i = \alpha_i(t_1) + A_i(t_2,\ldots,t_{g-1}) \quad i = 1, \ldots g
\leqno(1.1)$$\par

It is known that both the hypothesis and the conclusion of this theorem
effectively lead to geometric characterizations of Jacobian varieties among
all p.p.a.v.  For the hypothesis, this is work of Gunning and Welters.
See for example [W].  Thus, in what follows we will refer to the
hypothesis that $K(A)$ has a ``curve of flexes'' as {\it the Gunning-Welters
hypothesis}.  The subsequent paper [A-D] shows that this condition is
actually much {\it stronger} than necessary to characterize Jacobians.
Indeed, the existence of a jet of $V_{Y_p}$ of sufficiently high order
at one point is enough.  Naturally our proof will not make use of
any of these facts.  Instead we will connect the Gunning-Welters hypothesis
and the closely related solutions
$$u(x,y,t) = 2 {\partial^2\over\partial x^2} \log \theta_A(xU + yV + tW + z_0)
\leqno(1.2)$$
$ (U,V,W \in \bC^g) $  of the K.P. equation
$$(3/4)u_{yy} = (u_t - (1/4)(6 uu_x + u_{xxx}))_x \leqno(1.3)$$
to a differential-geometric criterion for the existence of parametrizations
of the form (1.1) for a hypersurface in $\bC^g$.  The proof of this
result occupies \S\S 2,3 of the present paper.

In one sense, our result is somewhat disappointing because it indicates that
the
approach to the Schottky problem via translation manifolds can be subsumed
in the other standard approaches. However, even though this is true, this
point of view does lead to one further characterization of Jacobian varieties
that may be of interest in its own right.  Namely, $(A, \Theta)$ is a Jacobian
if and only if its theta function satisfies a certain system of
fifth order PDE (found essentially by eliminating $U,V,W$ from the solution
of the K.P. equation given in (1.2)).  The equations are unfortunately
extremely complicated, but a procedure which describes how they may be
constructed may be found in [T].  We hope to return to these equations
in a future paper.

\beginsection 2. A Consequence of the Gunning-Welters Hypothesis and the
K.P. Equation\par
\bigskip
We begin by fixing some notation.
Let $A \cong \bC^g/(\bZ^g + \Omega\bZ^g)$ be a p.p.a.v. and let
$$\theta_A(z,\Omega) = \sum_{n\in \bZ^g} \exp \pi i[n^t\Omega n + 2n^tz]$$
be its Riemann theta function.  We will also make use of the related
theta functions with characteristics
$$\theta\left[\matrix{\eps\cr 0\cr}\right](2z,2\Omega) =
\sum_{n\in \bZ^g} \exp 2\pi i[(n+ \eps/2)^t\Omega (n + \eps/2)
                                + 2(n + \eps/2)^t z]$$
where $\eps$ is any vector of zeroes and ones in $\bZ^g$.  It is well-known
that thes functions form a basis for the linear system of second-order
theta functions on $A$.  We denote by $\vec{\theta_2}(z)$ the
vector-valued function
$$\vec{\theta_2}(z) = \left(\ldots,\theta\left[\matrix{\eps\cr
0\cr}\right](2z,2\Omega),\ldots\right)$$
The mapping $\psi_A$ in the Main Theorem is the projectivization of
$\vec{\theta_2}(z)$.

Following the work of Welters [W], we interpret the hypothesis of
the Main Theorem as follows.  Let $p \in \Gamma$.  The mapping
$Y_p \rightarrow A$ is induced by a local homomorphism
$$\eqalign{{\cal O}_{A,0} &\rightarrow \bC[t]/(t^3)\cr
                     f\quad &\rightarrow \sum_{i = 1}^2 \Delta_i(f)t^i\cr}$$
where the $\Delta_i$  are certain differential operators.  There are
corresponding
translation-invariant vector fields $D_1, D_2$ on $A$ such that formally
$$\exp\left(\sum_j D_j t^j \right) \equiv \sum_k \Delta_k t^k \pmod{t^3}$$
The Gunning-Welters hypothesis is that
$${\rm rank}\left( \vec{\theta_2}\ \Delta_1 \vec{\theta_2}\
                   \Delta_2\vec{\theta_2}\right) \le 2$$
at all points of ${1\over 2}(\Gamma - p)$.  We begin by noting

\proclaim Proposition 1.  ([W])  The artinian scheme $Y_p$ (defined by
$D_1$ and $D_2$) is equal to the second order jet of $V_{Y_p}$ at $p$.\par

This follows from the fact that for indecomposable p.p.a.v. $A$, the rank
of the $2^g \times (g(g+1)/2 + 1)$ matrix
$$\left(\vec{\theta_2}(0)\quad {\partial^2\over \partial z_i \partial z_j}
                               \vec{\theta_2}(0)\right)$$
is equal to $g(g+1)/2 + 1$.

Furthermore, the existence of a third-order jet of $\Gamma$ at $p$ implies
that there exist a vector field $D_3$ and a scalar $d \in \bC$ such that
$$0 = [D_1^4 - D_1D_3 + {3\over 4}D_2^2 + d]\vec{\theta_2}(0)\leqno(2.1)$$

We will write $D_1 = \sum_{i=1}^g U_i {\partial \over \partial z_i}$,
$D_2 = \sum_{i=1}^g V_i {\partial \over \partial z_i}$, and
$D_3 = \sum_{i=1}^g W_i {\partial \over \partial z_i}$.
By the Riemann quadratic theta formula, (2.1) can be rewritten as
$$0 = \theta_{xxxx}\theta - 4\theta_{xxx}\theta_x + 3\theta_{xx}^2
      + 4\theta_x\theta_t - 4\theta_{xt}\theta + 3\theta_{yy}\theta
      -3\theta_y^2 + 8d\theta^2,\leqno(2.2)$$
where $\theta = \theta(x,y,t) = \theta_A(xU + yV + tW + z_0)$.  By a
direct computation, this equation is seen to be equivalent to the equation
obtained by subsituting (1.2) into (1.3).

The key point for us will be that, under the assumption that
$\Gamma$ is one-dimensional, we can actually construct a two-dimensional
family of relations (2.2) for each $z_0 \in \bC^g$.  That is, there
is a 2-parameter family of triples $U,V,W$ which yield relations as in
(2.2).  This dimension count may be explained as follows.  We have one
parameter for the point $p \in \Gamma$, and a second for the choice
of $U \in T_p(\Gamma)$.

   The group of transformations
$$\eqalign{ U &\lra \lambda U\cr
            V &\lra \pm(\lambda^2 V + 2\alpha\lambda U)\cr
            W &\lra \lambda^3 W + 3\lambda^2\alpha + 3\lambda\alpha^2 U\cr
            d &\lra \lambda^4 d\cr}$$
leaves (2.1) invariant and acts on our set of triples $U,V,W$.  The
quotient projects to a curve $\Gamma^\prime$ in the projective space
$\bP^{g-1}$ with homogeneous coordinates $(U_1,\ldots,U_g)$.  In
other words, $\Gamma^\prime$ is the image of the projectivized Gauss
mapping on $\Gamma$.  (The relation between $\Gamma^\prime$ and $\Gamma$ is the
same as the relation between the canonical image of a curve $C$
and the $W_1$ subvariety of  $J(C)$.)

Now, let us specialize to the case $z_0 \in \Theta$ so that
$\theta(0,0,0) = \theta_A(z_0) = 0$.  Setting $x = y = t = 0$ in (2.2),
we have
$$0 = -4\theta_{xxx}\theta_x + 3\theta_{xx}^2 + 4\theta_x\theta_t -
       3\theta_y^2\leqno(2.3)$$
Furthermore, if we now add the extra condition that $U \ne 0$ is
{\it tangent} to $\Theta$ at $z_0$, so that $\theta_x(0,0,0) = 0$, then
(2.3) reduces to
$$0 = \theta_{xx}^2 - \theta_y^2$$
at $x = y = t = 0$, or
$$0 = \left(\sum_{i,j} \theta_{ij}(z_0)U_iU_j\right)^2 -
      \left(\sum_i \theta_i(z_0)V_i\right)^2,\leqno(2.4)$$
where we have written $\theta_i = {\partial \theta_A \over \partial z_i}$
and $\theta_{ij} = {\partial^2 \theta_A \over \partial z_i \partial z_j}$.

It is interesting to note that it is precisely this same special case
of the K.P. equation that was used by Beauville and Debarre in [B-D] to
link several different characterizations of Jacobians. We will use
(2.4) in a different way in what follows.  The condition that
$U \in T_{z_0}(\Theta)$ is that
$$\sum_i \theta_i(z_0)U_i = 0\leqno(2.5)$$
As $p$ varies on the curve $\Gamma \subset A$, the point with homogeneous
coordinates $(U_1,\ldots,U_g)$ varies on the curve
$\Gamma^\prime \subset \bP^{g - 1}$.  The equation (2.5) defines
the hyperplane $\bP T_{z_0}(\Theta) \subset \bP^{g-1}$, which will meet
the curve $\Gamma^\prime$ in a finite number of points
$q_1(z_0),\ldots,q_n(z_0)$.  If $z_0$ varies in a small open set on
$\Theta$, we can express the coordinates of any one of these intersection
points, say $q_1$, as analytic functions of any convenient set of local
coordinates on $\Theta$.  For simplicity, we will assume that
$U_1(z_0) \ne 0$ for all $z_0$ in our small open set on $\Theta$.

Since ${\rm dim}(\Gamma^\prime) = 1$, we may also assume (by renumbering
if necessary) that after we
divide by $U_1$ to obtain affine coordinates, each of
$\tau_2 = U_2/U_1,\ldots, \tau_g = U_g/U_1$
can be expressed in terms of $\tau_2$ on $\Gamma^\prime$.  We have
completed the preliminary constructions to prove

\proclaim Proposition 2.  Let $(A,\Theta)$ be an indecomposable p.p.a.v.
satisfying the Gunning-Welters hypothesis.  Then there exists an
analytic function $\tau_2$ on an open subset $X \subset \Theta$  and
analytic functions $\tau_3(\tau_2),\ldots,\tau_g(\tau_2)$, and
$\sigma_1(\tau_2),\ldots,\sigma_g(\tau_2)$ such that for all $z_0 \in X$,
\item{(i)}  $\sum_{i=1}^g \theta(z_0)\tau_i = 0$, (by convention,
here and in the following equation we put $\tau_1 = 1$),
\item{(ii)} $\sum_{i,j} \theta_{ij}(z_0)\tau_i\tau_j +
             \sum_i \theta_i(z_0)\sigma_i = 0$, and
\item{(iii)}  there is an analytic function $\lambda = \lambda(\tau_2)$ such
that $\sigma_i = \lambda\cdot {d\tau_i\over d\tau_2}$ for $i = 1,\ldots,g$.
\par

\noindent
{\bf Proof.}  (i) is the condition that $(U_1,\ldots,U_g) \in T_{z_0}(\Theta)$,
which holds by the argument given before
the statement of the Proposition.  From (2.4) and (i), we can obtain
relations such as (ii) for any $\sigma_i = \pm V_i/U_1^2 + \mu U_i$,
where $\mu$ is an arbitary function of $\tau_2$ (independent of $i$).  By
Proposition 1, at each point of $\Gamma$, $U$ and $V$ span the osculating
plane to $\Gamma$.  Hence, for the local parameter $\tau_2$ on $\Gamma$, we
have ${dU_i\over d\tau_2} = aV_i + bU_i$ for some functions $a,b$ of
$\tau_2$.  We can take
$$\eqalign{\sigma_i &=  {1\over a U_1}{d\tau_i\over d\tau_2}\cr
                    &= \left({1\over a U_1^2}{dU_i\over d\tau_2}
                    - {U_i\over a U_1^3}{dU_1\over d\tau_2}\right) \cr
                    &= {1\over U_1^2}V_i +
                       \left({b\over a} - {1\over a U_1^3}\right)U_i}$$
Then (ii) and (iii) follow from (2.4) and (i). $\triangle$
\bigskip
We remark that since (2.4) factors as
$$0 = \left(\sum_{i,j} \theta_{ij}(z_0)U_iU_j + \sum_i \theta_i(z_0)V_i\right)
 \cdot\left(\sum_{i,j} \theta_{ij}(z_0)U_iU_j - \sum_i \theta_i(z_0)V_i\right)
$$
there are usually at least {\it two} systems of functions
$\tau_i,\sigma_j$ satisfying the conditions of Proposition 2.  For if
$\tau_i,\sigma_j$ give one solution of the equations (i) and (ii), then
$-\tau_i,-\sigma_j$ give another.  This is a consequence of the {\it symmetry}
of the theta divisor.

\beginsection 3. Generalized Translation Manifolds\par
\bigskip
In this section, we will show that the conclusion of Proposition 2 of \S 2
implies that the theta-divisor $\Theta$ is a generalized translation
manifold (see the Main Theorem in \S 1 and [L]).  The
following proof is inspired by a similar discussion in [T].

\proclaim Proposition 3.  Let $H \subset \bC^g$ be an  analytic
hypersurface, defined by the equation $f(z_1,\ldots,z_g) = 0$.  Assume that
there exists an analytic function $\tau_2 = \tau_2(z_1,\ldots,z_g)$
on $H$, and analytic functions $\tau_1 = 1, \tau_i(\tau_2)$ for
$i = 3,\ldots, g$, and $\sigma_j(\tau_2)$ for $j = 1,\ldots,g$ satisfying
\item{(i)}  $\sum_i f_i(z)\tau_i = 0$,
\item{(ii)}  $\sum_{i,j} f_{ij}(z)\tau_i \tau_j +  \sum_i f_i(z) \sigma_i = 0$
and
\item{(iii)} there exists an analytic function $\lambda = \lambda(\tau_2)$
such that $\sigma_i = \lambda {d\tau_i \over d\tau_2}$
for $ i = 1, \ldots, g$.
\vskip 0truein
If $H$ is not {\it developable} (that is, if
the rank of the Gauss map on $H$ is generically $g - 1$), then
there exists an analytic parametrization
$$z_i = \alpha_i(t_1) + A_i(t_2,\ldots,t_{g-1})\qquad i = 1, \ldots,
g\leqno(3.1)$$
for $H$.  Conversely, the existence of a parametrization (3.1) for $H$
implies the existence of $\tau_i$ and $\sigma_j$ satisfying (i), (ii), and
(iii). \par

\noindent
{\bf Proof.}  The converse is easily seen by substituting (3.1) into the
equation of $H$ and differentiating twice with respect to $t_1$.  From
$$0 ={\partial\over\partial t_1} f(\alpha(t_1) + A(t_2,\ldots,t_{g-1}))$$
we obtain
$$0 = \sum_i f_i(z) \alpha_i^\prime(t_1)$$
for all $z \in H$.  Differentiating again,
$$0 =  {\partial^2\over\partial t_1^2} f(\alpha(t_1) + A(t_2,\ldots,t_{g-1}))
$$
which yields
$$0 = \sum_{i,j} f_{ij}(z) \alpha_i^\prime(t_1)\alpha_j^\prime(t_1)
      + \sum_i f_i(z) \alpha_i^{\prime\prime}(t_1)$$
Hence, by renumbering the coordinates if necessary, we can take $\tau_2 =
\alpha_2^\prime/\alpha_1^\prime$
and $\tau_i = \alpha_i^\prime/\alpha_1^\prime$ will all be functions
of $\tau_2$.  Then following the same idea as in the proof of Proposition
2, letting
$$\sigma_i = {1\over \alpha_1^\prime}{d\tau_i \over dt_1}
           = {1\over
(\alpha_1^\prime)^3}(\alpha_1^\prime\alpha_i^{\prime\prime}
              - \alpha_i^\prime\alpha_1^{\prime\prime})$$
we get a system of functions satisfying (i), (ii), and (iii).

For the direct implication, suppose that a system of functions $\tau_i,
\sigma_j$ exists satisfying (i), (ii), and (iii).  Consider
any one of the submanifolds $K \subset H$ defined by setting
$\tau_2 = t_0$ (constant). Let $z_0 = (z_{10},\ldots,z_{g0})$ be
an arbitrary point on $K$, and construct the integral curve $\alpha =
\alpha(z_1)$ of the system of ODE
$${dz_2 \over dz_1} = \tau_2, \ldots, {dz_g \over dz_1} = \tau_g$$
with initial condition $z_0$ in $\bC^g$.  By condition (i), along $\alpha$
$$\sum_i f_i(\alpha) {d\alpha_i \over dz_1} = \sum_i f_i(\alpha) \tau_i = 0
\leqno(3.2)$$
Hence $\alpha$ lies on $H$ (but not on $K$ or any of the other submanifolds
$\tau_2 = {\rm constant}$).

Differentiating (3.2) with respect to $z_1$, along $\alpha$ we have
$$0 = \sum_{i,j} f_{ij}(\alpha) \tau_i \tau_j
     + \sum_i f_i(\alpha) {d\tau_i \over d\tau_2}{d\tau_2\over dz_1}
\leqno(3.3)
$$
Subtract (3.3) from (ii) to obtain
$$0 = \sum_{i,j} f_i(\alpha)\left(\sigma_i - {d\tau_i \over
d\tau_2}{d\tau_2\over dz_1}\right)$$
By hypothesis (iii), we can write $\sigma_i = \lambda {d\tau_i \over d\tau_2}$
so this last equation becomes
$$0 = \left(\sum_i f_i(\alpha) {d\tau_i \over d\tau_2}\right)
      \left(\lambda - {d\tau_2 \over dz_1}\right)\leqno(3.4)$$
The first factor cannot be zero under our hypotheses since if it were,
$H$ would be developable.  Indeed from (3.3), we would have
$\sum_{i,j} f_{i,j}(\alpha)\tau_i \tau_j = 0$ and this implies that the
rank of the projectivized Gauss map is $\le g - 2$ everywhere on $X \subset H$,
as follows.  We express the submanifold $K$ locally as the intersection of
$f(z_1,\ldots,z_g) = 0$ and $z_g = h(z_1,\ldots, z_{g-1})$.  Then
for constant $\tau_2$, differentiating (i) with respect to $z_1,\ldots,z_{g-1}$
on $K$  we obtain
$$ \sum_i f_{ij}(z_0) \tau_i = -\sum_i f_{ig}(z_0) \tau_i {\partial h\over
\partial z_j} \quad j = 1, \ldots g -1 \leqno(3.5)$$
Since $\sum_j \tau_j \left[\sum_i f_{ij}(z_0)\tau_i\right] = 0$ by (3.3),
substituting from (3.5) we obtain
$$0 = \left[\sum_i f_{ig}(z_0)\tau_i\right]
    \left[\tau_g - \sum_{j=1}^{g-1} {\partial h\over\partial
z_j}\tau_j\right]$$
The second factor here is not zero since $\alpha$ is
not tangent to $K$.  Hence the first factor must be zero.  This together with
(3.5) shows that the projectivized Gauss map on $H$ has rank $\le g - 2$
at $z_0$.  Since $z_0$ was general on $K$, as $z_0$ and $\tau_2$ vary, we
see that $H$ would be developable if the first factor in (3.4) were zero.

Thus, from (3.4) $\lambda = {d\tau_2 \over dz_1}$, and the ODE defining the
curve $\alpha$ now take the form
$${dz_1\over d\tau_2} = {1\over \lambda},
 {dz_2 \over d\tau_2} = {\tau \over \lambda},\ldots,{dz_g\over d\tau_2} =
{\tau_g\over \lambda}$$
Integrating from $\tau_2 = t_0$ to $\tau_2 = t_1$ we obtain
a parametrization for $\alpha$ as follows.
$$z_i = z_{i0} + \int_{t_0}^{t_1} {\tau_i\over \lambda} d\tau_2\leqno(3.6)$$
$( i = 1, \ldots, g).$
The most important feature of these parametric equations is that the integrand
in
the second term is a function of $\tau_2$ {\it alone}.  As a consequence,
the integral curves $\alpha$ starting from all points $z_0 \in K$ are
{\it parallel translates} of one another in the ambient $\bC^g$.  This also
implies that the submanifolds $K$ are parallel translates of one another
in $\bC^g$.  Thus if we substitute any parametrization
$z_{i0} = A_i(t_2,\ldots,t_{g-1})$ for $K$ into (3.6), we obtain
a parametrization of the form (3.1).  The hypersurface generated in this
way must coincide with $H$, since as we noted above each integral
curve $\alpha$ is contained in $H$.  $\triangle$
\bigskip
Combining Propositions 2 and 3 completes the proof of the Main Theorem
from \S1, since the theta divisor of an indecomposable p.p.a.v. is
always non-developable.

In [L], we showed that a p.p.a.v. whose theta divisor has {\it two} distinct
parametrizations of the form (3.1) satisfying certain general position
hypotheses is the Jacobian of a non-hyperelliptic curve.  It seems highly
likely to us that these general position hypotheses are satisfied
automatically whenever a p.p.a.v. has two parametrizations of the form
(3.1), but we do not have a proof. The idea would be to analyze the
``degenerate'' generalized double translation manifolds where these
general position hypotheses do not hold, and show that they have
geometric properties incompatible with those of theta divisors, e.g. rulings,
etc.

The remark following Proposition 2
shows in addition that for theta divisors the existence of {\it one}
parametrization (3.1) usually implies the existence of a second.  It is
reasonable to conjecture that the only time this fails is for
hyperelliptic Jacobians, where the extra symmetry of the $W_1$ subvariety
causes the two generically distinct solutions of (2.5) to coincide.

We conclude by mentioning once again that sections 17-19 of
Tschebotarow's paper [T] contain a discussion of a systematic procedure
for eliminating the $\tau_i$, $\sigma_i$, $\lambda$ from the hypotheses
of Proposition 3.  The result is a simultaneous system of two fifth order PDE
on the defining equation of the hypersurface which are satisfied if and only
if the hypersurface is a generalized translation manifold.  Applied
to theta functions, these PDE should characterize the theta functions
of Jacobians.
\vfill\eject
{\bf References}
\bigskip
\frenchspacing
\item{[A-D]} Arbarello, E. and DeConcini, C. On a set of equations
characterizing Riemann matrices, {\sl Annals of Math.} {\bf 120} (1984),
119-141.
\smallskip
\item{[B-D]} Beauville, A. and Debarre, O. Une relation entre deux
approches du probleme de Schottky, {\sl Inv. Math.} {\bf 86} (1986),
195-207.
\smallskip
\item{[F]} Fay, J.  Schottky relations on ${1\over 2}(C - C)$, in
Proceedings of 1987 Bowdoin Summer Research Institute on Theta Functions,
{\sl Proc. of Symposia in Pure Math.} {\bf 49} (1989), 485-502.
\smallskip
\item{[L]} Little, J.  Translation manifolds and the Schottky problem, in
Proceedings of 1987 Bowdoin Summer Research Institute on Theta Functions,
{\sl Proc. of Symposia in Pure Math.} {\bf 49} (1989), 517-530.
\smallskip
\item{[M]}  Mumford, D. {\sl Curves and their Jacobians}  Ann Arbor: University
of Michigan Press, 1975.
\smallskip
\item{[T]}  Tschebotarow, N.  \"Uber Fl\"achen welche Imprimitivit\"atssysteme
in Bezug auf eine gegebene kontinuirliche Transformationsgruppe enthalten,
{\sl Receuil Math. (Sbornik)} {\bf 34} (1927), 149-204.
\smallskip
\item{[W]}  Welters, G. On flexes of the Kummer variety, {\sl Indagationes
Math.} {\bf 45} (1983), 501-520.
\bigskip
John B. Little
\smallskip
College of the Holy Cross
\smallskip
Worcester, MA 01610
\smallskip
U.S.A.
\end